\documentclass[a4paper,twoside,10pt]{amsart}
\usepackage{amsmath,amsfonts,amssymb,amsthm,mathtools}

\usepackage{color,graphicx}
\usepackage[a4paper]{geometry}
\usepackage[colorlinks=true,pdfstartview=FitH,pdfview=FitH]{hyperref}
\numberwithin{equation}{section}

\newcommand\mytop[2]{\genfrac{}{}{0pt}{}{#1}{#2}}

\def\expe{{\rm e}}

\def\iunit{{\rm i}}
\def\id{\,{\rm d}}
\def\bigO{{\mathcal O}}
\def\ifrac#1#2{\textstyle\frac{{#1}}{{#2}}\displaystyle}
\def\abs#1{\left|#1\right|}

\title{Exponentially-improved asymptotics and numerics for the (un)perturbed first Painlev\'e equation}
\author{Adri B. Olde Daalhuis}
\address{School of Mathematics and Maxwell Institute for Mathematical Sciences, The University of Edinburgh, Edinburgh EH9 3FD, United Kingdom}
\email{A.OldeDaalhuis@ed.ac.uk}
\urladdr{www.maths.ed.ac.uk/\string~adri} 

\dedicatory{Dedicated to Sir Michael V.\ Berry on the occasion of his 80$^{\text{th}}$ birthday.}

\begin{document}

\begin{abstract}
The solutions of the perturbed first Painlev\'e equation $y''=6y^2-x^\mu$, $\mu>-4$,
are uniquely determined by the free constant $C$ multiplying the exponentially
small terms in the complete large $x$ asymptotic expansions. 
Full details are given, including the nonlinear Stokes phenomenon, and the
computation of the relevant Stokes multipliers.
We derive asymptotic approximations, depending on $C$,
for the locations of the singularities that appear on the boundary of the sectors of
validity of these exponentially-improved asymptotic expansions. 
Several numerical examples illustrate the power of the approximations.
For the tri-tronqu\'ee solution of the unperturbed first Painlev\'e equation we give 
highly accurate numerics
for the values at the origin and the locations of the zeros and poles.
\end{abstract}

\maketitle

\section{Introduction}
The perturbed first Painlev\'e equation
\begin{equation}\label{ppy1}
    y''=6y^2-x^\mu,\qquad \mu\in{\mathbb R},~~\mu>-4,
\end{equation}
was discussed in \cite{Joshi2003}. In the case $\mu=0$ the solutions
are just Weierstrass $\wp$ functions, the case $\mu=1$ is the
unperturbed first Painlev\'e equation, and in the case $\mu=2$ two
solutions are $y(x)=\pm x/\sqrt6$, and no asymptotics is needed, although
the transseries solutions below are still valid. 

Before we start discussing the asymptotics of the solutions of \eqref{ppy1}
we first briefly discuss the possible singularities in the complex $x$ plane.
Obviously there could be a complicated branch-point at $x=0$. The other
singularities seem to be double poles, but a local analysis shows that
the local behaviour near such a point is of the form
\begin{align}
 y(x)&=\frac1{\left(x-x_j\right)^2}+\frac1{10} x_j^{\mu}\left(x-x_j\right)^2
 +\frac\mu{6} x_j^{\mu-1}\left(x-x_j\right)^3\nonumber\\ \label{polelog}
 &\qquad +\left(h_j-
 \frac{\mu(\mu-1)}{14} x_j^{\mu-2}\ln (x-x_j)\right)\left(x-x_j\right)^4\\
 &\qquad-\frac{\mu(\mu-1)(\mu-2)}{48} x_j^{\mu-3}\left(x-x_j\right)^5+\ldots,
 \nonumber
\end{align}
in which the only free constants are the location $x_j$ and the coefficient 
$h_j$. Note that the coefficient of $\left(x-x_j\right)^4$ does contain a
logarithm. This logarithm is absent in the known cases $\mu=0$ 
(Weierstrass $\wp$), $\mu=1$ (first Painlev\'e equation). 
In all other case these `poles' are actually branch-points.
Similar observations
have been made before. For example in \cite{Wittich53} and \cite{KC92} 
it is shown that for all of the solutions
of $y''(x)=6y(x)^2-f(x)$ to be single-valued about all movable singularities we need
$f''(x)=0$. In expansion \eqref{polelog} the next logarithm will appear in front of
$\left(x-x_j\right)^8$, and even higher powers of $\ln(x-x_j)$ will appear in the
tail of this expansion.

The dominant behaviours of the solutions that we will consider is $y_\pm(x)\sim \pm\sqrt{\frac{x^\mu}{6}}$.
In the next section we consider both behaviours, but in the remainder we will focus on $y_-(x)$. Note
that we have
\begin{equation}\label{plusminus}
   y_+(x)=\expe^{\frac{4\pi\iunit}{\mu+4}}y_-\left(x\expe^{\frac{2\pi\iunit}{\mu+4}}\right).
\end{equation}
Hence, our results for $y_-(x)$ can be translated to $y_+(x)$ via a rotation and a multiplication.

Rigorous results for the case $\mu=1$ are given in \cite{Joshi2001}. Some of
their techniques can also be applied in the case $\mu\not=1$. Proposition 2
of that paper gives us that for any sector of angle less than 
$\frac{4\pi}{\mu+4}$, there exist a solution $y(x)$ of \eqref{ppy1} 
such that $y(x)\sim-\sqrt{\frac{x^\mu}{6}}$ as $|x|\to\infty$ in this sector.
More importantly, their Theorem 3 shows us that there exist a unique solution
$y_-(x)$ of \eqref{ppy1} 
such that $y_-(x)\sim-\sqrt{\frac{x^\mu}{6}}$ as $|x|\to\infty$ in the sector
$|\arg x|<\frac{4\pi}{\mu+4}$. For this solution the constant beyond
all orders $C=0$.  It can also be labelled as being the
Borel-Laplace transform of its asymptotic expansion.

Once we have chosen the first term in the asymptotic expansion, the remaining
terms are fixed, and the free constant, say $C$, multiplies exponentially-small
terms. In \S\ref{Sectformal} we discuss the formal series solutions including
all the exponentially small terms in so-called transseries.
We will determine the sector of validity of these transseries.
The free constant $C$ will determine the location of the singularities 
that will appear on the boundary of the sector of validity.
At the end of \S\ref{Sectformal} we will give asymptotic approximations
for the locations of these singularities. In the numerical sections
of this paper it will be demonstrated that these approximations
are very good even for the singularities that are closest to the origin.

In \S\ref{Sectymin} we discuss the special solution $y_-(x)$, 
including the computation of its Stokes multipliers,
the level 1 hyperasymptotic approximation, which will include the
Stokes-smoothing of the Stokes phenomenon, and the location of
its singularities. In the numerical illustrations it seems to be the case
that this special solution has no singularities on the positive real
$x$-axis. This is known to be the case when $\mu=0$ and $\mu=1$.
As far as we can see the techniques of \cite{Joshi2001} can not
be used for $\mu\not=1$. 
The Borel transform of this formal series is discussed in \S\ref{SectBorel},
and it follows that this Borel-Laplace transform $y_-(x)$ is well defined for
$x>\tilde{\sigma}(\mu)$. For $\tilde{\sigma}(\mu)$ see Figure \ref{fig:sigma}.

The numerical tools are introduced in \S\ref{SectNumerics}. They are analytical continuation
via the Taylor-series method for analytic differential equations, and contour integral representations
for the locations of the singularities. These integrals are evaluated via the trapezoidal rule.
Note that once we have a reasonable guess for the location of a singularity, say $x_j$, and we can
evaluate $y(x)$ near that point, then \eqref{polelog} can also be used to obtain a much better
approximation for $x_j$.
These methods are used in \S\ref{SectExample1} for the cases $\mu=\frac{15}{7}$ and $\mu=4$.
Finally, in \S\ref{SectExample2} we illustrate the power of these simple methods by obtaining
approximations to a precision of $60$ significant digits for the unperturbed first Painlev\'e equation,
$\mu=1$, and in this way check some of the results in the recent literature.

The change of variable $z=\lambda x^{\frac\mu4+1}$ and
$y(x)=\sqrt{\frac{x^\mu}{6}} u(z)$, with $\lambda=\frac{8\cdot 6^{-\frac14}}{\mu+4}$,
will give us the differential equation
\begin{equation}\label{ppu1}
    u''(z)+2\nu\frac{u'(z)}{z}+\tfrac45 \nu\left(\tfrac65\nu-1\right)\frac{u(z)}{z^2}
    =\tfrac32\left(u^2(z)-1\right),
\end{equation}
in which we use the notation $\nu=\frac{5\mu}{2(\mu+4)}$. From an asymptotics point of view this
differential equation is slightly simpler than \eqref{ppy1}. 
Note that because we take $\mu>-4$ we will have $\nu<\frac52$.

\section{Formal series solutions}\label{Sectformal}
Differential equation \eqref{ppu1} has formal solutions of the form
\begin{equation}\label{level0}
    u_0(z)\sim \sum_{n=0}^\infty \frac{a_{n,0}}{z^n},
\end{equation}
with
\begin{equation}\label{coeff1}
    a_{0,0}=\pm1,\qquad a_{2n+1,0}=0,\qquad a_{2,0}=\tfrac4{15} \nu\left(\tfrac65\nu-1\right),\qquad
    a_{4,0}=2 \left(\tfrac2{15}\nu-1\right)\left(\tfrac35\nu-1\right)a_{0,0}a_{2,0},
\end{equation}
and the recurrence relation
\begin{equation}\label{coeff2}
    3a_{0,0}a_{n,0}=(n-2)\left(n-1-2\nu\right)a_{n-2,0}
    -\tfrac32\sum_{m=3}^{n-3}a_{m,0}a_{n-m,0},\qquad n=5,6,7,\ldots.
\end{equation}
This recurrence relation is consistent with $a_{2n+1,0}=0$.
We do give $a_{4,0}$ explicitly. It does not follow from the recurrence relation \eqref{coeff2}.

Our formal solution \eqref{level0} has no free constants. For the free constants we have to start
considering exponentially small perturbations for our solution, that is, solutions
of the form $u(z)=u_0(z)+Cu_1(z)$, in which $u_1(z)$ is exponentially small compared to $u_0(z)$.
However, our differential equations \eqref{ppy1} and \eqref{ppu1} are nonlinear and once we start
considering exponentially small terms we will immediately obtain terms that are double-, triple-, ...
exponentially small. Hence, we will end up with a transseries
\begin{equation}\label{trans}
    u(z)=\sum_{k=0}^\infty C^k u_k(z),
\end{equation}
in which the $u_k$ are solutions of the linear differential equations
\begin{equation}\label{ppuk}
    u_k''(z)+2\nu\frac{u_k'(z)}{z}+3\left(\frac{a_{2,0}}{z^2}-u_0(z)\right)
    u_k(z)=\tfrac32\sum_{\ell=1}^{k-1} u_\ell(z) u_{k-\ell}(z),\qquad k\geq1,
\end{equation}
with formal solutions
\begin{equation}\label{levelk}
    u_k(z)\sim \expe^{-k\sqrt{3a_{0,0}}z}\sum_{n=0}^\infty \frac{a_{n,k}}{z^{n+k\nu}}.
\end{equation}
In the case $k=1$ differential equation \eqref{ppuk} is linear and homogeneous, and hence,
there are no restrictions on $a_{0,1}$. We put that freedom in $C$ and fix $a_{0,1}=1$. For the
other coefficients we have
\begin{equation}\label{coeff3}
    2\sqrt{3a_{0,0}}n a_{n,1}=(n-1+\nu)(\nu-n)a_{n-1,1}
    +3\sum_{m=4}^{n+1}a_{m,0}a_{n-m+1,1},\qquad n=1,2,3,\ldots.
\end{equation}
The remaining coefficients are determined by
\begin{equation}\label{coeff4}
    (k^2-1)a_{0,0}a_{0,k}=\tfrac12\sum_{\ell=1}^{k-1}a_{0,\ell}a_{0,k-\ell},
    \qquad\Longrightarrow\qquad a_{0,k}=\frac{k}{\left(12a_{0,0}\right)^{k-1}},
    \qquad k=1,2,3,\ldots,
\end{equation}
and
\begin{align}
    &3(k^2-1)a_{0,0}a_{n,k}+2k\sqrt{3a_{0,0}}\left(n-1+(k-1)\nu\right)a_{n-1,k}\nonumber\\
    &\qquad+\left(n-2+k\nu\right)\left(n-1+(k-2)\nu\right)a_{n-2,k}
    -3\sum_{m=4}^{n}a_{m,0}a_{n-m,k}\nonumber\\
    &\qquad =\tfrac32\sum_{\ell=1}^{k-1}\sum_{m=0}^{n}a_{m,\ell}a_{n-m,k-\ell}.\label{coeff5}
\end{align}
Transseries expansion \eqref{trans} lives in the half-plane $\Re(\sqrt{3a_{0,0}}z)>0$. 
In this half-plane the terms decay exponentially, compare \eqref{levelk}. 
Below we will resum the transseries, and this is especially interesting on the
boundary of this sector.
In the opposite half-plane 
$\Re(\sqrt{3a_{0,0}}z)<0$ we have the transseries
\begin{equation}\label{trans2}
    u(z)=\sum_{k=0}^\infty C^k u_{-k}(z),
\end{equation}
in which
\begin{equation}\label{levelk2}
    u_{-k}(z)\sim \expe^{k\sqrt{3a_{0,0}}z}\sum_{n=0}^\infty \frac{\left(-1\right)^na_{n,k}}{z^{n+k\nu}}.
\end{equation}

When we combine \eqref{trans} with \eqref{levelk} we obtain a double
sum which can be resummed as
\begin{equation}\label{resum}
    u(z)\sim \sum_{n=0}^\infty \frac{G_n(X(z))}{z^{n}},
\end{equation}
in which $X(z)=C\expe^{-\sqrt{3a_{0,0}}z}z^{-\nu}$ and
\begin{equation}\label{resum2}
    G_n(X)= \sum_{k=0}^\infty a_{n,k}X^k.
\end{equation}
In the case of $n=0$ we can use \eqref{coeff4} to determine $G_0(X)$.
However, we can also substitute \eqref{resum} into \eqref{ppu1} and use
$X'(z)=-\left(\sqrt{3a_{0,0}}+\frac{\nu}{z}\right)X(z)$. This will give
us a power series expansion. The coefficient of $z^0$ can be evaluated as
\begin{equation}\label{G0eq}
   3a_{0,0}\left(X^2G_0''(X)+X G_0'(X)\right)=\ifrac32\left(G_0^2(X)-1\right),
\end{equation}
and the coefficient of $z^{-1}$ can be evaluated as
\begin{equation}\label{G1eq}
   3a_{0,0}\left(X^2G_1''(X)+X G_1'(X)\right)
   +2\nu\sqrt{3a_{0,0}}X^2 G_0''(X)=3G_0(X)G_1(X).
\end{equation}
Recall that $a_{0,0}^2=1$. We combine \eqref{G0eq} with the initial data
$G_0(0)=a_{0,0}$, $G_0'(0)=a_{0,1}=1$ and obtain
\begin{equation}\label{G0}
   G_0(X)=a_{0,0}+\frac{144X}{\left(X-12a_{0,0}\right)^2}.
\end{equation}
Similarly, we combine \eqref{G1eq} with the initial data
$G_1(0)=0$, $G_0'(0)=a_{1,1}=\nu(\nu-1)/(2\sqrt{3a_{0,0}})$ and obtain
\begin{equation}\label{G1}
   G_1(X)=\frac{\nu\sqrt{3a_{0,0}}X\left(288(1-\nu)-8(3\nu-19)a_{0,0}X
   +2X^2-\frac{a_{0,0}}{90}X^3\right)}{\left(X-12a_{0,0}\right)^3}.
\end{equation}
We already know from \eqref{polelog} that our solution can have double pole
type singularities. For a fixed $C$ we can use \eqref{G0} to obtain
a first approximation. The double poles should satisfy the approximation
$X(z)\approx12 a_{0,0}$. 

To obtain an extra term $X(z)\approx12 a_{0,0}+\frac{\alpha}{z}$ in this 
approximation we will look for double poles of $G_0(X)+z^{-1}G_1(X)$, that 
is, we determine constant $\alpha$ such that $G_0(X)+z^{-1}G_1(X)$ does not
have a triple pole at level $z^{-1}$. We expand
\begin{align}
   G_0(X)&=a_{0,0}+\frac{144X}{\left(X-12a_{0,0}
   -\frac\alpha{z}+\frac\alpha{z}\right)^2}\nonumber\\ \label{G0expand}
   &=a_{0,0}+
   \frac{144X}{\left(X-12a_{0,0}-\frac\alpha{z}\right)^2}
   -\frac{288X\alpha z^{-1}}{\left(X-12a_{0,0}-\frac\alpha{z}\right)^3}
   +\bigO{\left(z^{-2}\right)},
\end{align}
and determine $\alpha$ such that the triple pole in \eqref{G0expand} cancels
the triple pole in $z^{-1}G_1(X)$. The solution will be a function
of $z^{-1}$, but we are only interested in the constant part:
$\alpha=-\sqrt{3a_{0,0}}\nu\left(2\nu-\frac{124}{15}\right)$.
Hence, we expect double pole type singularities near solutions of
\begin{equation}\label{doublepoletype}
    C\expe^{-\sqrt{3a_{0,0}}z}z^{-\nu}=12 a_{0,0}
    -\sqrt{3a_{0,0}}\nu\left(2\nu-\ifrac{124}{15}\right)z^{-1}.
\end{equation}
The analysis above is similar to the one in \cite[\S6.6a]{CostinBook2009}.

\section{The case $y\sim - x^{\mu/2}/\sqrt{6}$}\label{Sectymin}
With the notation of the previous section we have $a_{0,0}=-1$,
and we take $\sqrt{3a_{0,0}}=\iunit\sqrt3$. Our starting point will be the
positive real axis and we consider the Borel-Laplace transform of the formal series
\begin{equation}\label{level0x}
    y_-(x)\sim \sqrt{\frac{x^\mu}6}\sum_{n=0}^\infty \frac{a_{n,0}\lambda^n}{x^{(\mu+4)n/4}},
\end{equation}
that is, the free constant $C=0$ when $x\to\infty$ along the positive real axis.
According to \eqref{levelk}  the `exponentially-small' terms are oscillatory on the positive 
real axis,  that is, the positive real axis is an anti-Stokes line. 
In the $z$-plane the imaginary axes will be active Stokes lines and the negative 
real axis will be the boundary for the sector of validity of asymptotic expansion 
\eqref{level0}. Hence, the sector
of validity is $\abs{\arg z}<\pi$, that is, $\abs{\arg x}<\frac{4\pi}{\mu+4}$.

The nonlinear Stokes phenomenon is the switching on of the exponentially small terms
when the imaginary axes are crossed, that is, in the transseries \eqref{trans} the constant
$C$ switches from $C=0$ to $C=K_\pm$ when we cross the positive/negative imaginary axis,
respectively.
The constants $K_\pm$ are the Stokes multipliers. The details
are very similar to the special case $\mu=1$ which is discussed in \cite{OD05b}. Hence, the
transseries expansions for this function are
\begin{equation}\label{ymintrans}
    y_-(x)\sim 
    \begin{cases}
    {\displaystyle\sqrt{\frac{x^\mu}6}\sum_{k=0}^\infty K_-^k u_k(z),} 
    &\frac{-4\pi}{\mu+4}<\arg x<\frac{-2\pi}{\mu+4},\\
     {\displaystyle\sqrt{\frac{x^\mu}6}u_0(z),} 
    &\frac{-2\pi}{\mu+4}<\arg x<\frac{2\pi}{\mu+4},\\
     {\displaystyle\sqrt{\frac{x^\mu}6}\sum_{k=0}^\infty K_+^k u_{-k}(z),}
    &\frac{2\pi}{\mu+4}<\arg x<\frac{4\pi}{\mu+4}.
    \end{cases}
\end{equation}
Near the boundaries of the sector of validity $\abs{\arg x}=\frac{4\pi}{\mu+4}$
the $u_k(z)$ are not exponential small anymore and it makes sense to resum the
transseries. Taking only the first term \eqref{levelk} and using \eqref{coeff4} we obtain
for $x$ near the boundary $\arg x=\frac{4\pi}{\mu+4}$ that
\begin{equation}\label{sing}
    y_-(x)\sim \sqrt{\frac{x^\mu}6}\sum_{k=0}^\infty 
    \frac{a_{0,k}K_+^k \expe^{\iunit k\sqrt3 z}}{z^{k\nu}}
    =\sqrt{\frac{x^\mu}6}\left(-1+
    \frac{K_+ \expe^{\iunit \sqrt3 z}z^{-\nu}}{\left(1+\frac1{12}K_+ \expe^{\iunit \sqrt3 z}z^{-\nu}\right)^2}.
    \right)
\end{equation}
Hence, the transseries contain information about singularities near the boundary
of the sector of validity. In this case we can see that we expect a double poles
near the solutions of $12+K_+ \expe^{\iunit \sqrt3 z}z^{-\nu}=0$. 
We know already from \eqref{polelog} that in the case $\mu\not=0,1$
these singularities are actually log singularities, but the dominant
term is a double pole. We will verify all of this in the
numerical sections below.

To determine the Stokes multipliers $K_\pm$ we can use the asymptotic formula
\begin{equation}\label{KpKm}
    a_{n,0}\sim \frac{K_+}{2\pi\iunit}\sum_{m=0}^\infty \left(-1\right)^m a_{m,1}
    \frac{\Gamma(n-m-\nu)}{\left(-\iunit\sqrt3\right)^{n-m-\nu}}
    -\frac{K_-}{2\pi\iunit}\sum_{m=0}^\infty a_{m,1}
    \frac{\Gamma(n-m-\nu)}{\left(\iunit\sqrt3\right)^{n-m-\nu}},
\end{equation}
as $n\to\infty$. Compare \cite[(4.3)]{OD05b}. Since $a_{2n+1,0}=0$ it follows that
\begin{equation}\label{KpKm2}
    K_+=\overline{K_-},
\end{equation}
and hence,
\begin{equation}\label{KpKm3}
    a_{2n,0}\sim \frac{-K_-}{\pi\iunit}\sum_{m=0}^\infty a_{m,1}
    \frac{\Gamma(2n-m-\nu)}{\left(\iunit\sqrt3\right)^{2n-m-\nu}},
\end{equation}
as $n\to\infty$. In this final result the optimal number of terms is $n$, and this formula
can be used to compute the Stokes multipliers numerically to any precision.

The details for the first hyperasymptotic re-expansion are very similar to the case $\mu=1$
discussed in \cite{OD05b}. The optimal number of terms of expansions \eqref{level0} and
\eqref{level0x} is $N$ such that $N-\sqrt3\abs{z}=\bigO(1)$ as $z\to\infty$. With this $N$
we have
\begin{equation}\label{Optlevel0}
    \sqrt{\frac6{x^\mu}}y_-(x)=\sum_{n=0}^{N-1}\frac{a_{n,0}}{z^n}+\bigO\left(\expe^{-\sqrt3\abs{z}}
    \abs{z}^{1/2}\right),
\end{equation}
as $z\to\infty$ in the sector $\abs{\arg z}<\frac12\pi$. The level 1 re-expansion will be
\begin{align}
    \sqrt{\frac6{x^\mu}}y_-(x)&=\sum_{n=0}^{2N-1}\frac{a_{n,0}}{z^n}
    +z^{1-2N}\frac{K_+}{2\pi\iunit}\sum_{n=0}^{N-1}\left(-1\right)^na_{n,1}
    F^{(1)}\left(z;\mytop{2N-n-\nu}{\iunit\sqrt3}\right)\nonumber\\ \label{Optlevel1}
    &\qquad -z^{1-2N}\frac{K_-}{2\pi\iunit}\sum_{n=0}^{N-1}a_{n,1}
    F^{(1)}\left(z;\mytop{2N-n-\nu}{-\iunit\sqrt3}\right)
    +\bigO\left(\expe^{-2\sqrt3\abs{z}}\abs{z}^{1}\right),
\end{align}
as $z\to\infty$ again in the sector $\abs{\arg z}<\frac12\pi$. 
Compare \cite[(5.3)]{OD05b}. The first hyperterminant function can
be expressed in terms of the incomplete gamma function 
$F^{(1)}\left(z;\mytop{N+1}{\sigma}\right)=-\expe^{\sigma z}\left(-z\right)^N\Gamma(N+1)
\Gamma(-N,\sigma z)$. It is the simplest function with a Stokes phenomenon. For more details
see \cite{OD98c}.

The level 1 expansion \eqref{Optlevel1} can be used to determine solution $y_-(x)$
uniquely. The order estimate in \eqref{Optlevel1} is double exponentially small.
Hence, the term $Cu_1(z)$ is clearly not present in the transseries expansion
for $y_-(x)$, that is, $C=0$.

The first array of poles in the lower half-plane are located near solutions of
\begin{equation}\label{polelowerminus}
K_-\frac{\expe^{-\iunit\sqrt3\lambda x^{(\mu+4)/4}}}{\lambda^\nu x^{\nu(\mu+4)/4}}
=-12-\frac{\iunit\sqrt3\nu\left(2\nu-\frac{124}{15}\right)}{\lambda x^{(\mu+4)/4}},
\end{equation}
and the first array of poles in the upper half-plane are located near solutions of
\begin{equation}\label{poleupperminus}
K_+\frac{\expe^{\iunit\sqrt3\lambda x^{(\mu+4)/4}}}{\lambda^\nu x^{\nu(\mu+4)/4}}
=-12+\frac{\iunit\sqrt3\nu\left(2\nu-\frac{124}{15}\right)}{\lambda x^{(\mu+4)/4}}.
\end{equation}
Note that the sign in front of the second term on the right-hand sides of
\eqref{polelowerminus} and \eqref{poleupperminus} are different.
This is a consequence of the $\left(-1\right)^n$, with $n=1$, in
\eqref{levelk2}.

\section{The numerics.}\label{SectNumerics}
To obtain very accurate numerical approximations we start with a large
$x_0$ on an anti-Stokes line and use an optimally truncated asymptotic
expansion. In the case of $y\sim -x^{\mu/2}/\sqrt6$ we will start
on the positive real $x$-axis and use \eqref{level0x} to determine
$y_-(x)$ and its derivative, and in the case $y\sim +x^{\mu/2}/\sqrt6$
we will start with a $x$ such that $\arg x=\frac{2\pi}{\mu+4}$. 
Once we have determined the function and its first derivative we can 
combine the original differential equation \eqref{ppy1} with the
Taylor-series method (see 
\cite[\href{http://dlmf.nist.gov/3.7.ii}{\S 3.7(ii)}]{NIST:DLMF})
and `walk' in the direction of the origin along the anti-Stokes line.
Thus initially we compute the first 2 Taylor coefficients in
\begin{equation}\label{Taylor}
    y(x)=\sum_{m=0}^\infty b_m\left(x-x_0\right)^m,
\end{equation}
via an optimally truncated asymptotic expansion, and the higher coefficients
via
\begin{equation}\label{Taylor2}
    (m+2)(m+1)b_{m+2}=6\sum_{\ell=0}^m b_\ell b_{m-\ell}-
    \left(-1\right)^mx_0^{\mu-m}\frac{\left(-\mu\right)_m}{m!},
    \qquad m=0,1,2,\ldots.
\end{equation}
We do control the step-size $step$ and the number of Taylor coefficients
that we use in \eqref{Taylor}. At each step we take $x_1=x_0+step$ and
compute $y(x_1)$ and $y'(x_1)$ via \eqref{Taylor} and take $x_0=x_1$ in
\eqref{Taylor2} to compute the higher coefficients.

To study the numerical stability of this process we can linearise
\eqref{ppy1} near $x_0$ and in that way we observe that the worst that can
happen is that the numerics is polluted with a little bit of a solution of 
\eqref{ppuk} ($k=1$). However, the solutions of that equation will be
oscillatory along the anti-Stokes line. Hence, the numerical integration should be
stable.

However, we are dealing with a nonlinear differential equation and do not
control the locations of the singularities. Note that in the previous sections
we did note that the origin can be a complicated branch-point and we did
make predictions of the locations of possible double `poles'.
Remarkably these predictions seem to be good even for small values of $x$.

In the case that $\mu=1$ we obtain from \eqref{polelog} that the residue at $x=x_j$ of
$-\frac12 xy'(x)/y(x)$ is $x_j$, and the reader can verify that the residue at $x=x_j$ of
$\frac1{56} y'(x)^3/y(x)$ is $h_j$. Hence, we can use loop integrals to evaluate
the position of the pole and the constant $h_j$. We do not know the exact location of the
poles, but we will need only reasonably good predictions, say $\tilde x_j$, 
which we do obtain from Pad\'e approximants, or from the solutions of \eqref{polelowerminus}
and \eqref{poleupperminus}.
Our loops will be circles $\abs{x-\tilde x_j}=r$ because we are going to use the trapezoidal rule,
and according to \cite{TW14} the right-hand side of
\begin{equation}\label{treff_eqn}
\frac{1}{2\pi\iunit}\oint_{|\tau|=r} F(\tau)\id \tau\approx\frac{1}{2M}
\sum_{m=0}^{2M-1}w_{m}F(w_{m}),\qquad\qquad\text{where}\quad w_{m}=r\expe^{\pi\iunit m/M},
\end{equation}
converges exponentially fast to the left-hand side as $M\to\infty$, 
as long as $(\tau-a)F(\tau)$ is analytic in a disc $|\tau|\leq \tilde r$, with $r<\tilde r$
and $|a|<r$. Once we know $y(x)$ and $y'(x)$ at, say, $x=\tilde x_j+w_0$ then we use 
\eqref{Taylor2} again to compute many Taylor coefficients, and use them in \eqref{Taylor} to 
evaluate $y(x)$ and $y'(x)$ at $x=\tilde x_j+w_1$. We can continue this process
to evaluate $y(x)$ and $y'(x)$ at all $x=\tilde x_j+w_m$.

In the case that $\mu\not=0,1$ this method to determine $x_j$ does not work, because $y(x)$
will have a logarithmic singularity at $x_j$. However, from \eqref{polelog} we obtain
the local expansion
\begin{equation}\label{rouche1}
    \frac{-xy'(x)}{2y(x)}=\frac{x_j}{x-x_j}+\ifrac{3}{14}\mu(\mu-1)x_j^{\mu-1}\left(x-x_j\right)^5
    \ln(x-x_j)+\cdots+reg(x-x_j),
\end{equation}
in which $reg(x-x_j)$ denotes a function that is analytic at $x=x_j$. Let $\tilde{x_j}$ be a
reasonable approximation for $x_j$ and let $r>0$ be small enough such that $x=x_j$ is the only
singularity contained in the disk $|x-\tilde{x_j}|\leq r$ then we can approximate the
integral
\begin{equation}\label{rouche2}
    \frac{1}{2\pi\iunit}\oint \frac{-xy'(x)}{2y(x)}\id x\approx x_j
    +\ifrac{1}{28}\mu(\mu-1)x_j^{\mu-1}\left(r+\tilde{x_j}-x_j\right)^6,
\end{equation}
where we integrate along the contour $x=\tilde{x_j}+r\expe^{2\pi\iunit\theta}$, $\theta\in[0,1]$.
Hence, the closer we are at $x_j$ the smaller the impact of
this logarithm. We will use \eqref{treff_eqn} with $F(\tau)=-\frac12 xy'(x)/y(x)$, in which
$x=\tau+\tilde{x_j}$, and decreasing values of $r$.

Note that we will ignore the second term on the right-hand side of \eqref{rouche2}.
However, the only unknown on the right-hand side of \eqref{rouche2} is $x_j$. Hence,
we can also evaluate the left-hand side of \eqref{rouche2} numerically, and use the full
approximation \eqref{rouche2} to compute $x_j$. This will result in slightly better
approximations.

\section{Example 1: $\mu=\frac{15}7$ and $\mu=4$}\label{SectExample1}
In the first example we take $\mu=\frac{15}7$, a non-integer. We will have
$\nu=\frac{75}{86}$. In this section we will aim to obtain approximations
to a precision of 10 significant digits. As a check we will do the same calculations
with larger starting points and considerably more steps and Taylor coefficients.
In this way we can check that the digits that we give below are actually correct.

\begin{figure}
    \centering
    \includegraphics[width=5.1truecm]{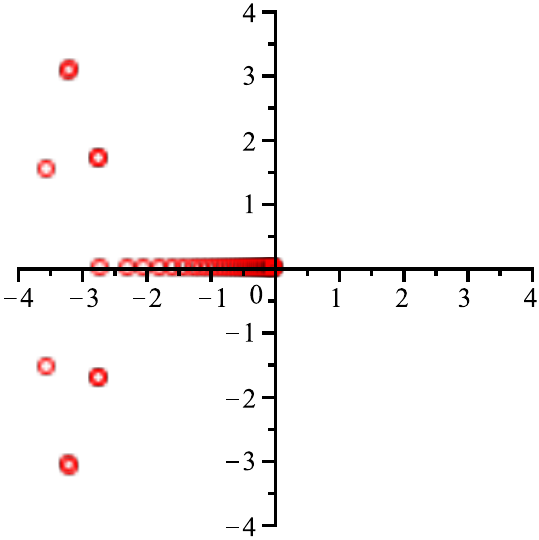}\qquad\qquad
    \includegraphics[width=5truecm]{figure1a.pdf}
    \caption{The poles of the Pad\'e approximants of $y_-(x)$ in the cases 
    $\mu=\frac{15}7$ (left) and $\mu=4$ (right).}
    \label{fig:Pade157}
\end{figure}

To determine the Stokes multipliers
we use \eqref{KpKm3} with $n=15$ and 15 terms on its right-hand side.
We obtain
\begin{equation}\label{Ex1K}
    K_-=0.07069725039+ 0.01439846034\iunit.
\end{equation}
For the remaining numerics in this section we use $x=6$ as our starting
point. The optimal number of terms in asymptotic approximation
\eqref{level0x} is approximately $11$. We obtain
\begin{equation}\label{Ex1b016}
    b_0=y_-(6)=-2.7837507946,\qquad b_1=y_-'(6)=-0.4971881751.
\end{equation}
The Taylor-series method is described in the previous section.
We will take $20$ Taylor coefficients in \eqref{Taylor} and `walk' in $100$
steps to $x=2$. We obtain
\begin{equation}\label{Ex1b012}
    y_-(2)=-0.8564979712,\qquad y_-'(2)=-0.4608802105.
\end{equation}
We compute the first $200$ Taylor coefficients at $x=2$
via \eqref{Taylor2} and use this Taylor series to compute a Pad\'e approximant
of order $[99,100]$ about the point $x=2$. In Figure \ref{fig:Pade157} (left) we can
see the distribution of the poles of this Pad\'e approximant. 
The accumulation of poles near the origin clearly indicates that the origin is a branch-point
(see \cite{Stahl97}), but we can also see that there are poles
at approximately $p_1\approx-2.75+1.7\iunit$ and at $p_2\approx-3.2+3.05\iunit$.

To obtain better numerical approximations for these poles we will `walk' along a straight line
from $x=2$ to $x=p_j+r$, with $r=\frac12$, and use the contour integral method described 
at the end of \S\ref{SectNumerics}.

In the case $j=1$, using $20$ Taylor coefficients we `walk' in
$1000$ steps to the first pole and obtain
\begin{equation}\label{Ex1b012}
    y_-(p_1+\ifrac12)=4.261757203 + 0.011948086\iunit,\qquad
    y_-'(p_1+\ifrac12)=-16.66671859 - 1.52622461\iunit.
\end{equation}
The size of these values indicates that we are close to the singularity.
The location of the pole is determined via \eqref{treff_eqn} in which we
take $F(\tau)=-\frac12 xy'(x)/y'(x)$, with $x=\tau+p_1$, 
$r=\frac12$ and $M=1000$. We will need $y(p_1+w_m)$. Again, we will use
the Taylor-series method with $20$ coefficients and step $w_m-w_{m-1}$.
We obtain the approximation
\begin{equation}\label{pole1}
    r=\ifrac12,\qquad p_1=-2.743854960 + 1.711273828\iunit.
\end{equation}
We repeat this process, starting at $x=2$, with this better guess for $p_1$ and
$r=\frac1{10}$
\begin{equation}\label{pole1b}
    r=\ifrac1{10},\qquad p_1=-2.740061378 + 1.709843142\iunit,
\end{equation}
and again
\begin{equation}\label{pole1c}
    r=\ifrac1{100},\qquad p_1=-2.740061121 + 1.709843110\iunit,
\end{equation}

At the end of \S\ref{Sectymin} we did mention that the poles should
approximately satisfy \eqref{poleupperminus}. This approximation was
constructed for the large poles. When we solve \eqref{poleupperminus}
for $x$ near $p_1$ we obtain the solution $x=-2.736 + 1.705\iunit$,
with relative error $0.002$.
Hence, even for the small poles we obtain reasonable approximations
via \eqref{poleupperminus}.

In a similar manner we can obtain a numerical approximation for the
second pole
\begin{align}\label{pole2}
    r=\ifrac1{2},\qquad   &p_2=-3.206143009 + 3.079481200\iunit,\nonumber\\
    r=\ifrac1{10},\qquad  &p_2=-3.200868582 + 3.074868336\iunit,\\
    r=\ifrac1{100},\qquad &p_2=-3.200868242 + 3.074868282\iunit,\nonumber
\end{align}
and when we solve \eqref{poleupperminus}
for $x$ near $p_2$ we obtain the solution $x=-3.199 + 3.074\iunit$,
with relative error $0.0004$.

In the introduction we do mention that the double pole expansion \eqref{polelog} can also be
used to obtain a very good approximation for the location of the pole. Say that we start  for the
second pole with the approximation originating from \eqref{poleupperminus}, that is 
$p_2\approx -3.199 + 3.074\iunit$ and with the Taylor series method, mentioned above, we evaluate
$y_-(p_2+\frac1{100})=6986.503356+1027.767205\iunit$. Then we can use the approximation
\begin{equation}\label{poleloglocation}
    y_-(x)\approx\frac1{\left(x-x_j\right)^2}+\frac1{10} x_j^{\mu}\left(x-x_j\right)^2
 +\frac\mu{6} x_j^{\mu-1}\left(x-x_j\right)^3,
\end{equation}
with $x=p_2+\frac1{100}$ and solve for $x_j$ near $p_2$. We obtain
$x_j=-3.200868241 + 3.074868282\iunit$. Note that compared with the final result in \eqref{pole2}
only the final digit is different.

We did mention above that in Figure \ref{fig:Pade157} (left) it is clearly visible
that the origin is a branch-point, but it is not obvious that the other poles are
actually also (weak) branch-points. For that reason we do include some details for the
case $\mu=4$, that is $\nu=\frac54$. 
In that case the origin is a regular point.
We compute the first $120$ Taylor coefficients at $x=0$
via \eqref{Taylor2} and use this Taylor series to compute a Pad\'e approximant
of order $[59,60]$ about the origin. In Figure \ref{fig:Pade157} (right) we can
see the distribution of the poles of this Pad\'e approximant. The poles start
to accumulate at $p_0$ and $p_1$, indicating that these are branch-points.
 With the same numerical steps
as described above we obtain that $K_-=0.5297382962 - 0.2194247868\iunit$, and that
there are `poles' at $p_0=-1.182001651$, $p_1=-0.895391503 + 2.352132859\iunit$,
and $p_2=-0.745388754 + 3.344311527\iunit$.

\section{Example 2: The unperturbed first Painlev\'e equation}\label{SectExample2}

In this section we will aim to obtain approximations
to a precision of 60 significant digits. The main reason for this is that we want to check
some of the results in the recent literature. The unperturbed first Painlev\'e equation
is the case $\mu=1$, $\nu=\frac12$. In this case the singularities in the complex plane
are double poles. Hence, they are not branch-points and the contour integral method
to determine the locations of the poles and zeros will be much more efficient.
The Stokes multiplier is known to be $K_-=-\frac{3^{1/4}}{\sqrt{5\pi}}(1+\iunit)$,
see \cite{Takei95}.
Taking $n=100$ in \eqref{KpKm3} and $100$ terms on its right-hand side we would
obtain an approximation for $K_-$ to a precision of 63 significant digits.

For the remaining numerics in this section we use $x=33$ as our starting
point. The optimal number of terms in asymptotic approximation
\eqref{level0x} is approximately $70$. We obtain
\begin{align}\label{Ex2b033}
    &b_0=y_-(33)=-2.345227006792405252263591282624246998603914831899264653960958,\\
    &b_1=y_-'(33)=-0.035532293810222842527936052573825449588186033237794348317154.\nonumber
\end{align}
We will take $40$ Taylor coefficients in \eqref{Taylor} and `walk' in $1000$
steps to the origin. We obtain
\begin{align}\label{Ex2b0origin}
    &y_-(0)=-0.187554308340494893838681757595444367707042203291560247736544,\\
    &y_-'(0)=-0.304905560261228856534104124988845544022671489625676976089364.\nonumber
\end{align}
Accurate values for this tri-tronqu\'ee solution at the origin are also given in
\cite{CD19}. They claim $64$ digit precision, but comparing their results with
\eqref{Ex2b0origin} we see that they did obtain $30$ digit precision.

\begin{figure}
    \centering
    \includegraphics[width=5truecm]{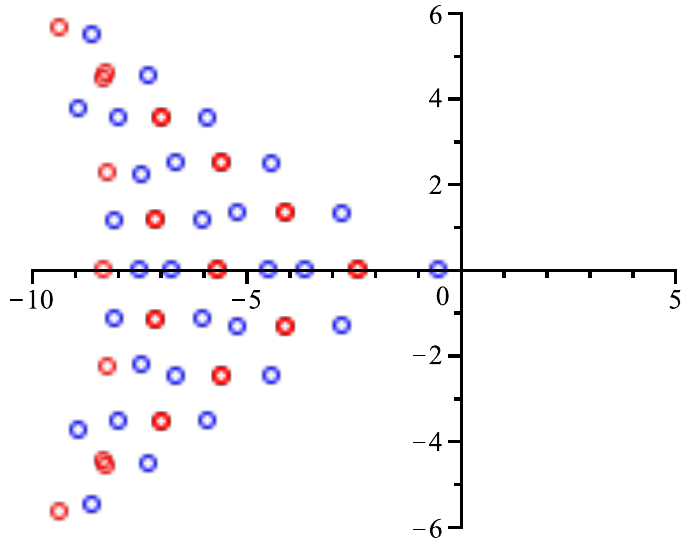}
    \caption{The zeros (blue) and poles (red) of the Pad\'e approximants of $y_-(x)$ in the case 
    $\mu=1$.}
    \label{fig:Pade1}
\end{figure}

We compute the first $100$ Taylor coefficients at $x=0$
via \eqref{Taylor2} and use this Taylor series to compute a Pad\'e approximant
of order $[49,50]$ about the point $x=0$. In Figure \ref{fig:Pade1} we can
see the distribution of the zeros and poles of this Pad\'e approximant. 

The location of the first zero, which is approximately at $z_1\approx-\frac12$, 
is determined via \eqref{treff_eqn} in which we
take $F(\tau)=\frac12 xy'(x)/y'(x)$, with $x=\tau-\frac12$, 
$r=\frac12$ and $M=60$. We will need $y(w_m-\frac12)$. Again, we will use
the Taylor-series method with $40$ coefficients and step $w_m-w_{m-1}$.
We obtain the approximation
\begin{align}\label{zero1}
    z_1&=-0.499912553551334521451561845356016137446077785951448892634807,\\ \nonumber
    y'(z_1)&=-0.468865514339593121531937054555736186201504711389139130341116,
\end{align}
in which we did obtain $y'(z_1)$ by walking in $10$ steps from the origin to $z_1$,
taking $40$ Taylor coefficients.
Note that with a relatively small $M=60$ we do already obtain more than
$60$ digits precision.

To approximate the first real pole we first walk to $x=-2$ in $300$ steps, taking
$40$ Taylor coefficients and obtain
\begin{align}\label{Ex2b0min2}
    &y_-(-2)=6.74868071988330557708652890818896015343487191457022993535053,\\
    &y_-'(-2)=-35.3975621098672136235306591226332932218101982620590565238107.\nonumber
\end{align}
The location of the first pole $p_1$ and the corresponding $h_1$ (see \eqref{polelog}), 
is determined via \eqref{treff_eqn} in which we
take $F(\tau)=-\frac12 xy'(x)/y'(x)$ and $\frac1{56} y'(x)^3/y(x)$, respectively, 
with $x=\tau-\frac52$, $r=\frac12$ and $M=200$. We obtain
\begin{align}\label{pole1mu1}
    p_1&=-2.38416876956881663929914585244876719041040881473785051267724,\\ \nonumber
    h_1&=0.0621357392261776408964901416400624601977407713738296636635327,
\end{align}
verifying the results in \cite{CD22}, except the final digits.
When we solve \eqref{poleupperminus}
for $x$ near $p_1$ we obtain the solution $x=-2.365+0.002\iunit$,
with relative error $0.008$.

Finally we compute also the location of the first complex pole. The details are the same
as above, the same number of steps, the same $M$, 
and the centre for the circle will be $-4.0+1.3\iunit$. The result is
\begin{align}\label{pole2mu1}
    p_2=&-4.07105552317228805392886956167452318934557741897847147742812\\ \nonumber
    &+1.33555121517567079951876062434077312552294901369825871527178\iunit,
\end{align}
and when we solve \eqref{poleupperminus}
for $x$ near $p_2$ we obtain the solution $x=-4.068+1.337\iunit$,
with relative error $0.0007$.

\section{The Borel transform on the positive real line}\label{SectBorel}
\begin{figure}
    \centering
    \includegraphics[width=6.5truecm]{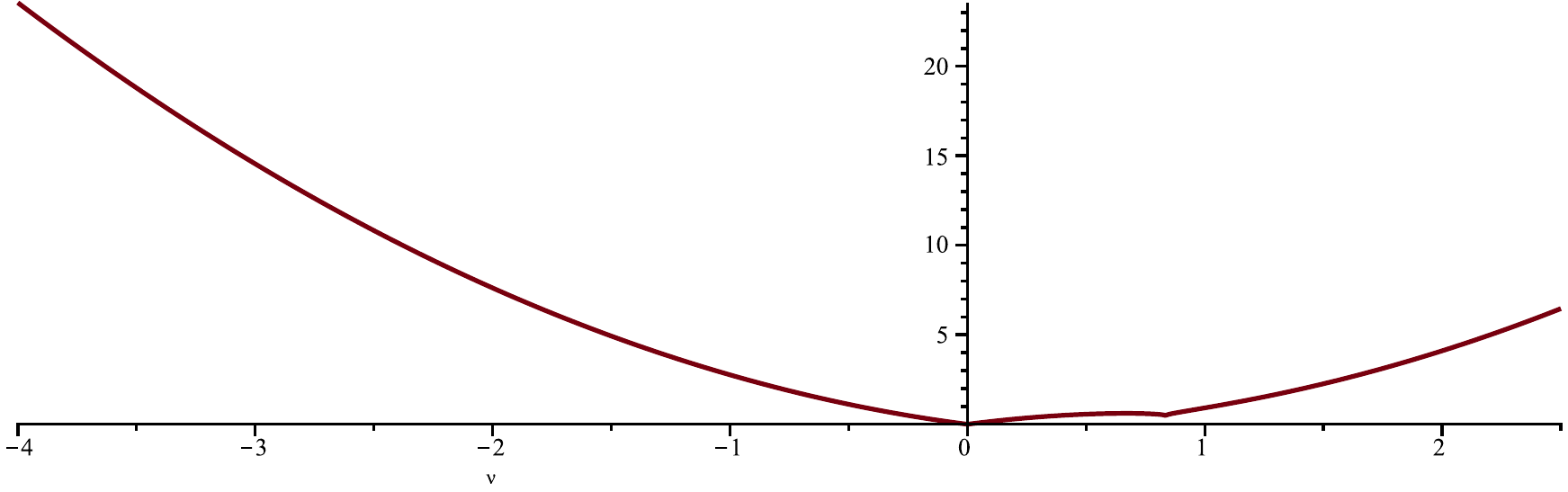}
    \qquad\includegraphics[width=6.5truecm]{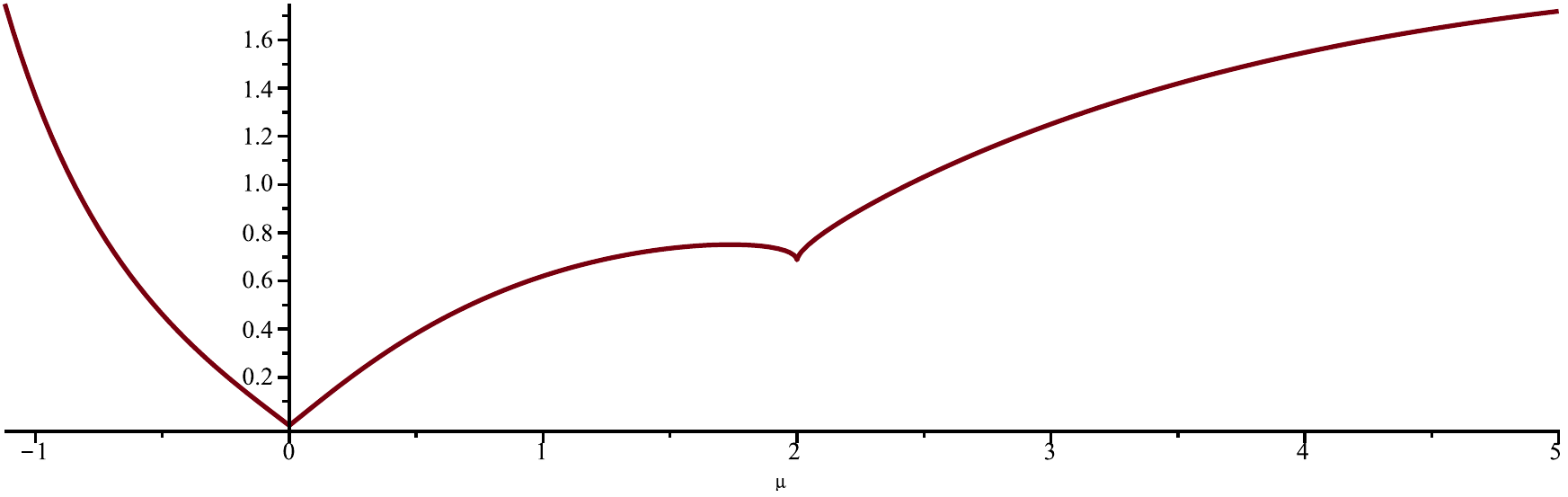}
    \caption{$\sigma$ as a function of $\nu$ (left) and
    $\tilde{\sigma}$ as a function of $\mu$ (right).}
    \label{fig:sigma}
\end{figure}

In this section we will study the Borel transform via
\begin{equation}\label{Borel}
    u(z)=-1+\int_0^\infty\expe^{-zt}b(t)\id t.
\end{equation}
When we start with differential equation \eqref{ppu1} and multiply each term
by $z^2$ then we obtain for the Borel transform $b(t)$ the 
differential-integral equation
\begin{equation}\label{integraleq1}
    \left(t^2+3\right)b''(t)+(2-\nu)2tb'(t)+2\left(\ifrac45\nu-1\right)
    \left(\ifrac35\nu-1\right)b(t)=\ifrac32
    \int_0^t b'(\tau)b'(t-\tau)\id \tau.
\end{equation}
It can be checked that
\begin{equation}\label{BorelExp}
    b(t)=\sum_{n=0}^\infty \frac{a_{n+1,0}}{n!}t^n,\qquad |t|<\sqrt3,
\end{equation}
with $a_{0,0}=-1$ and $a_{n+1,0}$ defined in \eqref{coeff1} and 
\eqref{coeff2}, is a solution of \eqref{integraleq1}. In this
section we will show that this solution is well defined and has a bound 
of the form $\left|b(t)\right|\leq c\expe^{\sigma(\nu) t}$, $t\geq0$.
Hence, our Borel-Laplace transform $u(z)$, defined in
\eqref{Borel}, is well-defined for $\Re(z)>\sigma(\nu)$. 
Our original perturbed first Painlev\'e equation is in terms of $x$
and $\mu$, and we obtain that the corresponding solution is well defined
for $x>\tilde{\sigma}(\mu)=
\left(\frac{6^{1/4}}{8}(\mu+4)\sigma\right)^{1/(1+\mu/4)}$.
For $\sigma(\nu)$ and $\tilde{\sigma}(\mu)$ see Figure \ref{fig:sigma}.

To obtain a more convenient integral equation
we use \eqref{ppu1} directly in \eqref{Borel} and obtain
\begin{equation}\label{integraleq2}
    \left(t^2+3\right)b(t)=2\nu\int_0^t \tau b(\tau)\id \tau
    +3a_{2,0}t- 3a_{2,0}\int_0^t (t-\tau) b(\tau)\id \tau
    +\ifrac32\int_0^t b(\tau)b(t-\tau)\id \tau.
\end{equation}
Dividing both sides by $t^2+3$ we have $b(t)={\mathcal T}b(t)$ with
\begin{equation}\label{integraleq}
    {\mathcal T}b(t)=\frac{2\nu}{t^2+3}\int_0^t \tau b(\tau)\id \tau
    +\frac{3a_{2,0}t}{t^2+3}
    - \frac{3a_{2,0}}{t^2+3}\int_0^t (t-\tau) b(\tau)\id \tau
    +\frac{\frac32}{t^2+3}\int_0^t b(\tau)b(t-\tau)\id \tau.
\end{equation}
We are going to show that this is a contraction mapping.
Let $c$ and $\sigma$ be positive constants and define the norm
\begin{equation}\label{norm}
    \Vert h\Vert=\inf \left\{ M~\vert~|h(t)|\leq Mc \expe^{\sigma t}
    ~{\rm for~all}~t\geq0\right\}.
\end{equation}
Denote by ${\mathcal B}_\sigma$ the complex vector space of analytic function $h(t)$ on $[0,\infty)$ 
such that $\Vert h\Vert$ is bounded. Equipped with this norm, ${\mathcal B}_\sigma$ becomes a Banach space.

For the terms on the right-hand side of \eqref{integraleq} we have in 
the case $t\geq0$ that
\begin{equation}\label{B0}
    \left|\frac{2\nu}{t^2+3}\int_0^t \tau h(\tau)\id \tau\right|
    \leq\frac{2|\nu|c\Vert h\Vert}{t^2+3}\int_0^t\tau\expe^{\sigma\tau}
    \id\tau
    \leq\frac{2|\nu|c\Vert h\Vert t}{t^2+3}\int_0^t\expe^{\sigma\tau}
    \id\tau\leq \frac{|\nu|\Vert h\Vert}{\sigma\sqrt{3}}c\expe^{\sigma t},
\end{equation}
\begin{equation}\label{B1}
    \left|\frac{3a_{2,0}t}{t^2+3}\right|
    \leq |a_{2,0}|\left(t\expe^{-\sigma t}\right)\expe^{\sigma t}
    \leq \frac{|a_{2,0}|}{c\sigma}c\expe^{\sigma t},
\end{equation}
\begin{equation}\label{B2}
    \left|\frac{3a_{2,0}}{t^2+3}\int_0^t (t-\tau) h(\tau)\id \tau\right|
    \leq\frac{3|a_{2,0}|t}{t^2+3}\int_0^t |h(\tau)|\id \tau
    \leq \frac{\sqrt3|a_{2,0}|\Vert h\Vert}{2\sigma}c\expe^{\sigma t},
\end{equation}
\begin{equation}\label{B3}
    \left|\frac{\frac32}{t^2+3}\int_0^t h_1(\tau)h_2(t-\tau)\id\tau\right|
    \leq\frac{\frac32 t\Vert h_1\Vert\Vert h_2\Vert}{t^2+3}c^2\expe^{\sigma t}
    \leq \frac{\sqrt3\Vert h_1\Vert \Vert h_2\Vert c}{4}c\expe^{\sigma t},
\end{equation}
in which we have used several times $\frac{t}{t^2+3}\leq\frac{\sqrt3}{6}$.

Combining the inequalities above we obtain
\begin{equation}\label{B5}
    \Vert{\mathcal T}h\Vert\leq
    \frac{|a_{2,0}|}{c\sigma}
    +\ifrac{\sqrt3}{4}c\Vert h\Vert^2 +
    \frac{|\nu|+\ifrac32|a_{2,0}|}{\sigma\sqrt3}\Vert h\Vert,
\end{equation}
and using the convolution identity $h_1*h_1-h_2*h_2=(h_1+h_2)*(h_1-h_2)$ 
we have
\begin{equation}\label{B6}
    \Vert{\mathcal T}h_1-{\mathcal T}h_2\Vert\leq
    \left(\ifrac{\sqrt3}{4}c\Vert h_1+h_2\Vert +
    \frac{|\nu|+\ifrac32|a_{2,0}|}{\sigma\sqrt3}\right)
    \Vert h_1-h_2\Vert.
\end{equation}
Recall that $a_{2,0}=\tfrac4{15} \nu\left(\tfrac65\nu-1\right)$. 
It is now possible to choose $c$ and $\sigma$ such that when
$\Vert h\Vert\leq1$
we will have from \eqref{B5} that $\Vert{\mathcal T}h\Vert\leq1$,
and taking $\Vert h_j\Vert\leq1$ we will have in \eqref{B6} that the
multiplier of $\Vert h_1-h_2\Vert$ will be less than $1$. 
Hence, we want to find a pair $c,~\sigma$ such that both
\begin{equation}
    \frac{|a_{2,0}|}{c\sigma}
    +\ifrac{\sqrt3}{4}c +
    \frac{|\nu|+\ifrac32|a_{2,0}|}{\sigma\sqrt3}\leq1,
    \qquad
    \ifrac{\sqrt3}{2}c +
    \frac{|\nu|+\ifrac32|a_{2,0}|}{\sigma\sqrt3}\leq1.
\end{equation}
Once we have such a pair we have shown that $b(t)={\mathcal T}b(t)$
has a unique solution with the bound $\left|b(t)\right|\leq c\expe^{\sigma t}$
for all $t\geq0$. 
We want $\sigma$ as small as possible,
but $\sigma\sim\frac4{25}\left(\sqrt{3}+\frac2{c}\right)\nu^2$ 
as $\nu\to-\infty$. A reasonable choice seems to be $c=\frac7{10}$.
It is also possible to use the optimal $c=\sqrt{\frac3{\alpha^2}+\frac4{\alpha}}-\frac{\sqrt{3}}{\alpha}$,
with $\alpha=\left|\frac{\nu}{a_{2,0}}\right|+\frac32$.
The corresponding $\sigma$ as a function of $\nu$ is displayed in
Figure \ref{fig:sigma}.

\section*{Acknowledgement}
The author wants to thank Nalini Joshi for stimulating discussions regarding the main topics 
of this paper, and thanks the Isaac Newton Institute for Mathematical Sciences for support 
during the program ‘Applicable resurgent asymptotics: towards a universal theory’ 
supported by EPSRC grant no.~EP/R014604/1.
The authors’ research was supported by a research grant 60NANB20D126 from the National Institute of Standards and Technology.

\bibliographystyle{siam}
\bibliography{temp}

\end{document}